\input amstex

\documentstyle{amsppt}
\magnification=1200
\hcorrection{.25in}
\advance\vsize-.75in

\document
\topmatter

\title{Tiling the integers with translates of one finite set}
\endtitle

\author{Ethan M. Coven and Aaron Meyerowitz}
\endauthor

\address{Ethan M. Coven,
         Deparment of Mathematics,
         Wesleyan University,
         Middletown, CT  06459-0128}
\endaddress

\email{ecoven\@wesleyan.edu}
\endemail

\address{Aaron  Meyerowitz,
         Deparment of Mathematics,
         Florida Atlantic  University,
         Boca Raton, FL  33431-0991}
\endaddress

\email{meyerowi\@fau.edu}
\endemail

\thanks
{Part of this work was done while the first author was a member
of the Mathematical Sciences Research Institute (MSRI), where research
is supported in part by NSF grant DMS-9701755.  The authors thank J.~
Propp for putting them in touch with each other.  The first author
thanks J.~Jungman and M.~ Keane for helpful conversations.}
\endthanks

\abstract
A set {\it tiles the integers\/} if and only if the integers
can be written as a disjoint union of translates of that set.  We
consider the problem of finding necessary and sufficient conditions
for a finite set to tile the integers.  For sets of prime power size,
it was solved by D.~Newman [J.  Number Theory {\bf 9} (1977),
107--111].  We solve it for sets of size having at most two prime
factors.  The conditions are always sufficient, but it is unknown
whether they are necessary for all finite sets.
\endabstract


\define \poly{polynomial}
\define \cyclo{cyclotomic}
\define \nneg{nonnegative}
\define \lcm{\operatorname{lcm}}

\define \coeffs{coefficients}
\define \coeff{coefficient}
\define \st{such that}

\endtopmatter

\head Introduction
\endhead

Let $A$ be a finite set of integers.  $A$ {\it tiles the integers\/}
if and only if the integers can be written as a disjoint union of
translates of~$A$, equivalently, there is a set $C$ such that every
integer can be expressed uniquely $a+c$ with $a \in A$ and $c \in C$.
In symbols, $A \oplus C=\Bbb Z$.  In this case $A$ is called a {\it
tile\/}, $A \oplus C=\Bbb Z$ a {\it tiling\/}, and $C$ the {\it
translation set\/}. For a survey of such tilings, see R.~Tijdeman
\cite{Tij}. For connections with group theory and functional analysis,
see \cite{Haj} and \cite{L-W}.

We consider the problem of finding necessary and sufficient conditions
for a finite set to tile the integers.  For sets of prime power size
(cardinality, denoted $\#$), it was solved by D.~Newman \cite{New}.
Newman remarked that ``even for so simple a case as size six we do not
know the answer.'' We find necessary and sufficient conditions for~$A$
to tile the integers when $\#A$ has at most two prime factors.

There is no loss of generality in restricting attention to translates
of a finite set $A$ of {\it nonnegative\/} integers.
Then $A(x) = \sum_{a\in A}x^a$ is a \poly\ \st\  $\# A=A(1)$.
Let $S_A$ be the set of prime powers~$s$ such that
the $s$-th \cyclo\ \poly\ $\Phi_s(x)$ divides~$A(x)$. Consider the
following conditions on $A(x)$.

\roster
 \widestnumber\item{(T2)}
  \item"{(T1)}" $A(1) = \prod_{s\in S_A} \Phi_s(1)$.
  \item"{(T2)}" If $s_1,\dots,s_m\in S_A$ are powers of distinct
                primes, then $\Phi_{s_1\cdots s_m}(x)$ divides~$A(x)$.
\endroster

\proclaim{Theorem A}   
If $A(x)$ satisfies  (T1) and (T2), then $A$ tiles the integers.
\endproclaim

\proclaim{Theorem B1} 
If $A$ tiles the integers, then $A(x)$ satisfies (T1).
\endproclaim

\proclaim{Theorem B2} 
If $A$ tiles the integers and $\# A$ has at most two prime factors,
then $A(x)$ satisfies ~(T2).  
\endproclaim

\proclaim{Corollary} 
If $\# A$ has at most two prime factors, then $A$ tiles the integers if
and only if $A(x)$ satisfies (T1) and (T2).
\endproclaim

It is unknown whether the sufficient conditions (T1) and (T2) are
necessary for any finite set to tile the integers. (T1) is necessary
but not sufficient (see the example after Theorem ~B1 in
Section~2). However, if $\# A$ is a prime power, then (T2) follows
from~(T1), so in this case (T1) {\it is\/} necessary and sufficient. An
examination of Newman's proof \cite{New,~Theorem~1} essentially yields
this result.

Our proof of Theorem~B2 provides a structure theory for finite sets
~$A$ \st\ $A$ tiles the integers and $\#A$ has at most two prime
factors.  We sketch this in Section~4.

If $A$ is a finite set which tiles the integers, then
$ \bigcup_{a \in A} [a,a+~1)$ tiles the reals.
J.~ Lagarias and Y.~ Wang \cite{L-W} proved a structure theorem for
closed subsets~$T$  of the reals with finite Lebesgue measure and
boundary of measure zero \st\ the reals can be written as a countable
union of measure-disjoint translates of~$T$. It describes such sets
in terms of finite sets which tile the integers.

\head 1. Preliminaries \endhead

For $A$ and $B$ sets or multisets of integers, we denote the
multiset $\{a+b: a \in~A, b \in B\}$ by $A+B$. We write $A \oplus B$
when every element can be expressed uniquely   $a+b$.
For $k$ an integer, we write $kA$ for $\{ka: a \in A\}$,
we call $\{k\} \oplus A$ a {\it translate \/} of $A$,
and when $k$ is a factor of every $a \in A$, we write ~$A/k$ for
$\{a/k: a \in A\}$.

For $s \ge 1$, the $s$-th {\it cyclotomic polynomial\/}
$\Phi_s(x)$ is defined recursively by $x^s - 1 =  \prod\Phi_t(x)$,
where the product is taken over all factors~$t$ of~$s$.  
The factors of~$s$ are positive and include both $1$ and~$s$.
\proclaim{Lemma 1.1}  Let $p$ be prime.  Then
\roster
  \item $\Phi_s(x)$ is the minimal \poly\ of any primitive $s$-th root
        of  unity.
  \item $1+x+\dots+x^{s-1} = \prod\Phi_t(x)$, where
        the product is taken over all factors $t > 1$ of~$s$.
  \item $\Phi_p(x) = 1 + x + \dots +
        x^{p-1}$ and $\Phi_{p^{\alpha +1}}(x)=\Phi_p(x^{p^\alpha})$.
  \item $\Phi_s(1) = \cases
                  0&\text{if $s=1$}\\
                  q&\text{if $s$ is a power of a prime $q$}\\
                  1&\text{otherwise}.
                  \endcases$
  \item $\Phi_s(x^p)= \cases
                 \Phi_{ps}(x)&\text{if $p$ is a factor of~$s$}\\
                 \Phi_s(x)\Phi_{ps}(x)&\text{if $p$ is not a factor
                  of~$s$}.
                 \endcases$
  \item If $s$ and $t$ are relatively prime, then
        $\Phi_s(x^{t}) = \prod\Phi_{rs}(x)$,
        where the product is taken over all factors $r$ of~$t$.
          
  \item If $\bar A(x)$ is a polynomial and $A(x)=\bar A(x^p)$, then
        $ \{t:\Phi_t(x) \text{ divides } A(x)\}=$\newline
        $\{s':\Phi_s(x) \text{ divides } \bar A(x)\} \cup
        \{ps:\Phi_s(x) \text{ divides } \bar A(x)\}$, where $s'=ps$ or
        $s$ according as $p$ is or is not a factor of $s$.
\endroster
\endproclaim

\demo{Proof}
(1) is a  standard fact.
(2) and (3) follow from the definition,
(4) from~(2) and~(3),
and (5) from (1) because the roots of $\Phi_s(x^p)$ are
$e^{2\pi ik/ps}$ for $k$ relatively prime to $s$.
Repeated application of~(5) yields~(6).
For (7), let $\omega=e^{2 \pi i/t}$. Then $\omega^p$ is  a primitive
$s$@-th root of unity for some $s$ and, from~(5),
$t \in \{s',ps\}$. $\Phi_t(x)$ divides $\bar A(x^p)$ if and only 
$\bar A(\omega^p)=0$ if and only if
$\Phi_s(x)$ divides $\bar A(x)$.
\qed\enddemo

A set $C$ of integers is {\it periodic\/} if and only if $C \oplus
\{n\}=C$ for some $n \ge 1$.
Then $C$ is a union of congruence classes modulo $n$ and $C=B\oplus n
\Bbb Z$, where $B$ is any set consisting of one representative from
each class.  If $A \oplus C = \Bbb Z$ is a tiling and $C$ is periodic,
the smallest such~$n$ is called the {\it period\/} of the tiling.  Note
that $n=(\#A)(\#B)$ and $A \oplus B$ is a complete set of residues
modulo ~$n$. Conversely, if $A \oplus B$ is a complete set of residues
modulo~$n$, then $A \oplus (B \oplus n\Bbb Z) = \Bbb Z$
is a tiling of period~$n$ or less, as are 
$B \oplus (A \oplus n\Bbb Z) = \Bbb Z$ and 
$A' \oplus C = \Bbb Z$ for any
$A' \equiv A \pmod n$.

The following basic result is due to G\.~Haj\'os \cite{Haj} and
N\.~deBruijn \cite{deB-1}, then C\.~Swenson \cite {Swe}, then Newman
\cite{New}.

\proclaim{Lemma 1.2} 
Every tiling by translates of a finite set is periodic, i.e., if $A$
is  finite  and $A \oplus C\ = \Bbb Z $, then there is a finite set~$B$
such that $C = B \oplus n \Bbb Z$, where $n = (\#A)(\#B)$. \qed
\endproclaim

\remark{Remark}
Newman's proof shows that the period of any tiling by~$A$ is bounded
by $2^{\max (A) - \min (A)}$.  The tiling $\{j\} \oplus \Bbb Z = \Bbb
Z$ has period~1. The tiling $A \oplus C=\Bbb Z$, where $A=\{j\} \oplus
\{0,k\}$ and $C=\{0,1,\dots,k-1\}\oplus 2k\Bbb Z$, has period $2k$.
We know of no other tilings whose period is as large as $2\left( \max
(A) - \min (A)\right)$.  See the remarks following Lemma 2.1.
\endremark

The collection of all finite multisets of \nneg\ integers is in
one-to-one correspondence with the set of all \poly s with \nneg\
integer \coeffs.  The correspondence is $$A \longleftrightarrow A(x) =
\sum_{a \in A}m_ax^a,$$ where $m_a$ is the multiplicity of $a$ as an
element of $A$. If $B$ is another such multiset and $k \ge 1$, then
the polynomial corresponding to~$A+B$ is~$A(x)B(x)$,
to~$A \cup B$ is~$A(x) + B(x)$, and to~$kA$ is~$A(x^k)$.  Using this
language we get

\proclaim{Lemma 1.3}
Let $n$ be an integer and let
$A$ and $B$  be finite multisets of \nneg\ integers  with
corresponding \poly s $A(x)$ and ~$B(x)$.
Then the following statements are
equivalent. Each forces $A$ and $B$ to be sets
\st\ $(\#A)(\#B)=A(1)B(1)=n$.
\roster
 \item $ A \oplus (B \oplus n\Bbb Z)=\Bbb Z$ is a tiling.
 \item $A \oplus B$ is a complete set of residues modulo~$n$.
 \item $A(x)B(x) \equiv 1 + x +  \dots + x^{n - 1}
         \pmod{x^n-1}$.
 \item    $n=A(1)B(1)$ and for every factor $t > 1$
       of~$n$, the \cyclo\ \poly\
       $\Phi_t(x)$ is a divisor of~$A(x)$ or~$B(x)$. \qed
\endroster
\endproclaim

There is no loss is restricting attention to conditions for a finite
set of {\it nonnegative\/} integers to tile the integers.
We can further restrict to finite sets whose minimal element is~$0$
and to translation sets which contain~$0$, although we will
not always do so. For if $A'$ and ~$C'$ are translations of~$A$ and
~$C$, then $A \oplus C=\Bbb Z$ if and only if $A' \oplus C'=\Bbb Z$.

Recall that~(T1) and~(T2) concern the set $S_A$ of prime powers $s$ such
that the cyclotomic polynomial $\Phi_s(x)$ divides $A(x)$.
When $A$ and a translate $A'$ are finite sets of nonnegative integers,
$A(x)$ and $A'(x)$ are divisible by the same cyclotomic polynomials, so
\roster
 \item"{$\bullet$}"$A$ tiles the integers if and only if $A'$ tiles
                   the integers.
 \item"{$\bullet$}"$A(x)$ satisfies (T1)
                   if and only if $A'(x)$ satisfies (T1).
 \item"{$\bullet$}" $A(x)$ satisfies (T2)
                    if and only if $A'(x)$ satisfies (T2).
\endroster

The next lemma allow us to further restrict attention to finite sets of
integers with greatest common divisor~$1$.

\proclaim{Lemma 1.4} Let $k > 1$ and let
$A=k\bar A$ be a finite set of nonnegative integers.
\roster
 \item  $A$ tiles the integers if and only if $\bar A$ tiles the
        integers.
 \item  If $p$ is prime, then 
        $S_{p\bar A} = \{p^{\alpha+1} : p^\alpha \in S_{\bar A}\}
        \cup \{q^\beta \in S_{\bar A}\: q {\text{ prime\/}}, q \ne p\}$.
 \item $A(x)$ satisfies (T1) if and only if $\bar A(x)$ satisfies (T1).
 \item $A(x)$ satisfies (T2) if and only if $\bar A(x)$ satisfies (T2).
\endroster
\endproclaim

\demo{Proof}
For one direction of (1), let 
$\bar A \oplus C=\Bbb Z$. 
Then $k\bar A \oplus kC = k\Bbb Z$ and hence
$A \oplus (\{0,1,\dots,k-1\} \oplus kC) = \Bbb Z$.
For the other, let $k\bar A \oplus D =\Bbb Z$.
Then $k\bar A \oplus D_0 =k\Bbb Z$, where $D_0 = \{d \in D :
d\equiv 0 \pmod k\}$, and hence $\bar A \oplus D_0/k = \Bbb Z$.
(2) follows from Lemma~ 1.1(7).

It suffices to prove~(3) and (4) when $k$ is prime, say $k=p$. (3)
follows from~(2) and Lemma~1.1(4)
since $\# A = \#\bar A$.
It remains to prove (4).  Let $s' =ps$ or~$s$
according as $p$ is or is not a factor of~$s$.
Let  $s_1,\dots,s_m$ be powers of distinct primes and $s=s_1\cdots
s_m$. Then  $s_1',\dots, s_m'$ are powers of distinct primes and
$s'=s_1'\cdots s_m'$. From~(2), every $s_i \in S_{\bar A}$ if
and only if every $s_i' \in S_A=S_{p \bar A}$. From~1.1(7), $\Phi_s(x)$
divides $\bar A(x)$ if and only if $\Phi_{s'}(x)$ divides $A(x)$.
Putting all this together yields~(4).
\qed
\enddemo

\remark{Remark} It follows from (2) that
$B$ is not contained in $p\Bbb Z$ when $\Phi_p(x)$ divides $B(x)$.
\endremark
\medskip

Lemma~1.4 deals with $A \subset k\Bbb Z$.  The related situation that
$A \oplus C=\Bbb Z$ is a tiling with  $C \subseteq k\Bbb Z$
leads to an important construction.  We defer it to Lemma~2.5.

\head 2. Tiling results  \endhead

\proclaim{Theorem A} Let $A$ be a finite set of \nneg\ integers with
corresponding \poly\ $A(x) = \sum_{a\in A}x^a$ and let $S_A$ be the set
of prime powers~$s$ such that the \cyclo\ \poly\ $\Phi_s(x)$
divides~$A(x)$.  If
\roster
 \widestnumber\item{(T2)}
  \item"{(T1)}" $A(1) = \prod_{s\in S_A} \Phi_s(1)$.
  \item"{(T2)}" If $s_1,\dots,s_m \in S_A$ are powers of distinct
        primes, then  $\Phi_{s_1\cdots s_m}(x)$ divides~$A(x)$,
\endroster
then $A$ tiles the integers.
\endproclaim

\demo{Proof}
We construct a set~$B$ such that condition~(4) of Lemma~1.3 is
satisfied.  Define $B(x) = \prod\Phi_s(x^{t(s)})$, where the product
is taken over all prime power factors~$s$ of~$\lcm(S_A)$ which are not
in ~$S_A$ and  $t(s)$ is the largest factor of~$\lcm(S_A)$ relatively
prime to~$s$.  Since every such~$s$ is a prime power, $B(x)$ has \nneg\
\coeffs. Since 1.3(4) will be shown to hold, these \coeffs\ are all $0$
and~$1$.

Let $s > 1$ be a factor of~$A(1)B(1)$ and write $s=s_1 \cdots s_m$ as
a product of powers of distinct primes.  If every $s_i \in S_A$, then
by~(T2), $\Phi_s(x)$ divides~$A(x)$.  Suppose then that some $s_i
\notin S_A$.  Then $\Phi_{s_i}(x^{t(s_i)})$ divides~$B(x)$, $r=s/s_i$
is a factor of $t(s_i)$ and, by Lemma~1.1(6) (with $s=s_i$ and
$t=t(s_i)$), $\Phi_{rs_i}(x)$ divides $\Phi_{s_i}(x^{t(s_i)})$.
Thus  $\Phi_s(x)$ divides~$B(x)$ since $rs_i=s$.
\qed
\enddemo

\remark{Remarks} The set $B$ constructed in the proof depends only on
$S=S_A$ and not on~$A$.  Defining $C_S=B\oplus \lcm (S)\Bbb Z$, 
$A \oplus C_S=\Bbb Z$ for {\it all \/} $A$ with $S_A=S$ which
satisfy~(T1) and~(T2). Then $C_S \subseteq p\Bbb Z$ for every prime $p
\in S$, since $p$ is a factor of~$n$ and
every divisor $\Phi_s(x^{t(s)})$ of $B(x)$ is a polynomial in $x^p$.
For either $t(s)$ is a multiple of $p$, or $s=p^{\alpha+1}$  
with $\alpha \ge 1$  and
$\Phi_s(x^{t(s)})=\Phi_p\left(x^{t(s)p^\alpha}\right)$,
so every  divisor $\Phi_s(x^{t(s)})$ of $B(x)$ is a
polynomial in $x^p$.
\endremark

\proclaim{Theorem B1} 
Let $A$ be a finite set of \nneg\ integers with corresponding
\poly\ $A(x) = \sum_{a\in A}x^a$ and let $S_A$ be the set of prime
powers~$s$ such that the \cyclo\ \poly\ $\Phi_s(x)$ divides~$A(x)$. If
$A$ tiles the integers, then
 \roster 
  \item"{(T1)}" $A(1) = \prod_{s \in S_A} \Phi_s(1)$.
 \endroster
\endproclaim

\remark{Remark} (T1) is
not sufficient for~$A$ to tile
the integers. $A = \{0,1,2,4,5,6\}$ does not tile the
integers, but $A(x) = \Phi_3(x)\Phi_8(x)$ satisfies~(T1).
\endremark

\medskip

Theorem B1 follows from Lemma 2.1(1) below.

\proclaim{Lemma 2.1} 
Let $A(x)$ and $B(x)$ be polynomials with  \coeffs\ $0$ and ~$1$, 
$n=A(1)B(1)$, and $R$   the set of prime power factors of~$n$.
If $\Phi_t(x)$ divides $A(x)$ or~$B(x)$
for every factor $t > 1$ of~$n$, then
 \roster
  \item $A(1)=\prod_{s\in S_A}\Phi_s(1)$ and 
        $B(1)=\prod_{s\in S_B}\Phi_s(1)$.
  \item $S_A$ and $S_B$ are disjoint sets whose union is $R$.
 \endroster
\endproclaim

\demo{Proof}
For every factor $t > 1$ of~$n$, $\Phi_t(x)$ divides $A(x)$ or~$B(x)$,
so $R \subseteq S_A \cup S_B$. Clearly $A(1) \ge \prod_{s \in
S_A}\Phi_s(1)$ and $B(1) \ge \prod_{s \in S_B}\Phi_s(1)$. Thus
$$A(1)B(1)\ge \prod_{s\in S_A}\Phi_s(1) \prod_{s\in
  S_B}\Phi_s(1)\ge \prod_{t\in R}\Phi_t(1) = n,
$$
the equality by Lemma~1.1(4). Hence all the inequalities and
containments above are actually equalities, and $S_A$ is disjoint from
~$S_B$.
\qed \enddemo

\remark{Remarks} If a tiling $A \oplus C =\Bbb Z$ has period $n$ and
$C=B \oplus n\Bbb Z$, then $n=\lcm(S_A \cup S_B)$, so the period of any
tiling by $A$ is a multiple of $\lcm(S_A)$. A particular tiling by~$A$
may have period larger than $\lcm(S_A)$, however when $A(x)$ satisfies
(T1) and (T2), the tiling
$A \oplus \left(B \oplus (\#A)(\#B)\Bbb Z\right) = \Bbb Z$
constructed in the proof of Theorem A has period~$\lcm(S_A)$.  In all
cases known to the authors both $A(x)$ and $B(x)$ satisfy (T1) and (T2).

We leave it to the interested reader to show that for any set $A$
of \nneg\ integers,

\roster
 \item"{$\bullet$}"  $ \lcm(S_A) \le\frac{p}{p-1}\left(\max (A)-\min
                     (A)\right)$, where $p$ is the smallest prime
                     factor of~$\#A$.
 \item"{$\bullet$}" The inequality is strict except when 
       $A=\{j\} \oplus p^{\alpha}\{0,1,\dots,p-~1\}.$
\endroster

We show in Lemma 2.3 that there is always a tiling whose
period is a product of powers of the prime factors of~$\#A$.  
\endremark

\proclaim{Theorem B2} Let $A$ be a finite set of \nneg\ integers with
corresponding \poly\ $A(x) = \sum_{a\in A}x^a$  such that $\#A$ has at
most two prime factors and let $S_A$ be the set of prime powers~$s$
such that the \cyclo\ \poly\ $\Phi_s(x)$ divides~$A(x)$. If $A$ tiles
the integers, then
 \roster
  \widestnumber\item{(T2)}
   \item"{(T2)}" If $s_1,\dots,s_m \in S_A$ are powers of distinct
    primes, then $\Phi_{s_1\cdots s_m}(x)$ divides~$A(x)$.
\endroster
\endproclaim

The following result is crucial to our proof of Theorem~B2. We give an
alternate proof of it in Section ~3.

\proclaim{Lemma 2.2} \cite{Tij, Theorem~1} 
Suppose  that $A$ is finite, $0 \in A \cap C$,
and $A \oplus C = \Bbb Z$. If $r$ and $\# A$ are relatively prime, then
$rA \oplus C = \Bbb Z$. \qed
\endproclaim

\remark{Remark}
Translating $A$ or $C$ does not affect the conclusion.
Thus the condition $0 \in A \cap C$ is not needed.
\endremark

\proclaim{Lemma 2.3}  
If a finite set $A$ tiles the integers, then there is a tiling by~ $A$
whose period is a product of powers of the prime factors of $\#A$.
\endproclaim

\demo{Proof} 
If $A\oplus C= \Bbb Z$ is a tiling of period $n$ and $r > 1$ is a
factor of $n$ relatively prime to $\#A$, then by Lemma~2.2, 
$ rA \oplus C=\Bbb Z$. Therefore $ rA \oplus C_0 =r\Bbb Z$, where 
$C_0 = \{c \in C : c \equiv 0 \pmod{r}\}$, and hence 
$ A \oplus C_0/r = \Bbb Z$ is a tiling of period~$n/r$. 
\qed
\enddemo

The following result is essentially Theorem~4 of \cite{San}.  We prove a
more general result which implies it in Section~3.

\proclaim{Lemma 2.4} \cite{San} 
Let $A \oplus C=\Bbb Z$ be a tiling of period~$n$ such that $A$ is 
finite, $0 \in A \cap C$, and $n$ has one or two prime factors. Then
there is a prime factor $p$ of~$n$ such that either $A \subset p\Bbb Z$
or $C \subseteq p\Bbb Z$. \qed \endproclaim

Sands' result is stated in the terms of direct sum decompositions of
finite cyclic groups, but it is easy to translate it into the
terminology of this paper.

\proclaim{Lemma 2.5} Suppose $A \oplus C=\Bbb Z$, where $A$ is a finite
set of nonnegative integers, $k>1$, and $C \subseteq k\Bbb Z$.  For
$i = 0,1,\dots, k-1$, let $A_i=\{a \in A : a \equiv i \pmod k\}$,
$a_i=\min(A_i)$, and $\bar A_i=\{a-a_i : a \in A_i\}/k$.  Then

 \roster
  \item $A(x)=x^{a_0}\bar A_0(x^k)+x^{a_1}\bar A_1(x^k)+\dots
        +x^{a_{k-1}}\bar A_{k-1}(x^k)$.
  \item Every $\bar A_i \oplus C/k=\Bbb Z$.
  \item The elements of $A$ are equally distributed modulo $k$ --- every
        $\#\bar A_i=(\#A)/k.$
  \item $S_{\bar A_0}=S_{\bar A_1}=\cdots=S_{\bar A_{k-1}}$.
  \item When $k$ is prime, $S_A=\{k\} \cup S_{k\bar A_0}$ and
        if every $\bar A_i(x)$ satisfies (T2), then $A(x)$
        satisfies~(T2).
 \endroster  
\endproclaim

\demo{Proof} (1) is clear.  (2)~follows from
$A_i \oplus C=\{i\} \oplus k\Bbb Z=\{a_i\} \oplus k\Bbb Z$.
To prove~(3), note that the  translation set $C/k$
has some period~$n$, so there is a set $\bar B$ such that
$\bar A_i \oplus (\bar B \oplus n\Bbb Z)=\Bbb Z$
and every $\bar A_i \oplus \bar B$
is a complete set of residues modulo~$n$.
Thus the $\#\bar A_i$ are equal, so (3) holds.
(4)~also follows since by Lemma~2.1, every $S_{\bar A_i}$ is
the complement of $S_{\bar B}$ in 
the set of
prime power factors of~$n$.

To prove (5), write $p$ in place of ~$k$. From Lemma~1.4(2), 
$S_{p\bar A_i} = \{s' : s \in S_{\bar A_i}\}$, where $s' = ps$ or~$s$
according as $p$ is or is not a factor of~$s$. The polynomial
corresponding to $p\bar A_i$ is $\bar A_i(x^p)$, so from~(1) and~(4),
$S_{p\bar A_0} \supseteq S_A$.  Also $p \in S_A$, since if
$\Phi_p(\omega) = 0$, then $\omega^p = 1$, $\omega^{a_i}=\omega^i$, and
$A(\omega) = \sum_{i=0}^{p-1}\omega^i \bar A_i(1) = (\# A/k)
\sum_{i=0}^{p-1}\omega^i = 0$, 
the next-to-last equality by~(3). We have thus shown that $S_A
\supseteq \{p\} \cup S_{p\bar A_0}$.  Since $A_0$ and~$A$ tile the
integers, $A_0(x)$ and $A(x)$ satisfy~(T1) and
$S_A = \{p\} \cup S_{p\bar A_0}$.

Now assume that every $\bar A_i(x)$ satisfies~(T2).
Condition (T2) for $A(x)$  is: if $s_1,\dots ,s_m \in S_{\bar A_0}$
are powers of distinct primes, then $\Phi_{s_1'\cdots s_m'}(x)$ divides
~$A(x)$ and $\Phi_{ps_1\cdots s_m}(x)$ divides~$A(x)$.
By~(T2), $\Phi_{s_1\cdots s_m}(x)$ divides every~$\bar A_i(x)$.
Hence by Lemma~1.1(7), $\Phi_{s_1'\cdots s_m'}(x)$ and
$\Phi_{ps_1\cdots s_m}(x)$ divide all the ~$\bar A_i(x^p)$, so they
divide $A(x)$ as well.
\qed \enddemo

\proclaim{Corollary} 
If $A$ is a finite
set of integers
 and $C \subseteq k\Bbb Z$, then $A \oplus C=\Bbb Z$
if and only if  
$A=\bigcup_{i=0}^{k-1}\left(\{a_i\} \oplus k \bar A_i \right)$for
some complete set $\{a_0,a_1,\dots,a_{k-1}\}$ of residues modulo~$k$,
and $k$ sets $\bar A_i$, each of which satisfies
$\min(A_i)=0$ and tiles the integers with translation set $C/k$.
\qed
\endproclaim

The decomposition is unique. We can have $\gcd(A)=1$ although this
may not be true for the ~$\bar A_i$. If the $\bar A_i$ are
equal, then the union is a direct sum,
$A~=~\{a_0,a_1 \dots,a_{k-1}\} \oplus k\bar A_0$.
For some simple choices of translation set~$C$, every tile has this
form.

\demo{Proof of Theorem B2} 
From Lemma 1.4 and the comments before it
there is no loss of generality in assuming that $\gcd(A)=1$
and $0\in A$.

By Lemma 2.3 there is a tiling $A \oplus C=\Bbb Z$ whose
period~$n$ is a product of powers of the prime factors of
$\#A$. We complete the proof by induction on ~$n$.
If  $n=1$, then $A=\{0\}$ and $A(x)\equiv 1$ satisfies (T2) vaccuously.
If $n >1$, then by Lemma~2.4, there is a prime factor $p$ of $n$ such
that $C \subseteq p\Bbb Z$. Then by Lemma~2.5,
$A(x)=x^{a_0}\bar A_0(x^p)+x^{a_1}\bar
A_1(x^p)+\dots +x^{a_{p-1}}\bar A_{p-1}(x^p)$
and every $\bar A_i\oplus C/p=\Bbb Z$ is a tiling of period~$n/p$. By
the inductive hypothesis, every $\bar A_i(x)$ satisfies (T2), so by
Lemma~ 2.5(5), $A(x)$ satisfies~(T2). 
\qed \enddemo

Every set known to the authors, regardless of size, which tiles the
integers satisfies the tiling conditions (T1) and~(T2).  However, our
proof of Theorem~B2 cannot be extended to sets whose size has more than
two prime factors because Lemma~2.4 need not hold.
For $m$  a positive integer with more than two prime factors,
a very general construction due to S.~Szab\'o \cite{Sza}
gives sets~$A$ \st\ $\#A=m$, $\min(A)=0$, $\gcd(A)=1$, and $A$ tiles
the integers, yet the members of $A$ are {\it not \/} equally
distributed modulo $k$ for any $k>1$. Hence, from Lemma~2.5(3), every
set~$C$ \st\ $0 \in C$ and $A \oplus C=\Bbb Z$ satisfies $\gcd(C)= 1$.
All these sets $A$ satisfy~(T1) and~(T2).

These examples also show that both Tijdeman's
conjecture \cite{Tij, p.~266} --- if $A \oplus C=\Bbb Z$, 
$0\in A \cap C$, and $\gcd(A)=1$, then $C \subseteq p\Bbb Z$ for some
prime factor of $\#A$ --- and the weaker conjecture
--- if $A$ tiles the integers, $\min(A)=0$ and $\gcd(A)=1$,
then there is {\it some \/} translation set
of the desired type --- are false without further conditions.
Tijdeman's conjecture would have implied an inductive characterization
of all tilings $A \oplus C= \Bbb Z$. The weaker conjecture would have
implied  an inductive characterization of the finite sets which tile
the integers. We established the weaker conjecture in Lemma~2.4 for
those $A$ \st\ $\#A$ has one or two prime factors.
We show how to use it in Section~4. Tijdeman \cite{Tij, Theorem 3}
proved his conjecture when $\# A$ is a prime power. We do not know
whether it holds when $\# A$ has exactly two prime factors.

\head 3. Alternate proofs of Tijdeman's and Sands' Theorems
\endhead

Tijdeman's  Theorem (Lemma 2.2) follows from Lemma~1.3 and

\proclaim{Lemma 3.1}
Let $A$ and $B$ be finite sets of \nneg\ integers with
corresponding \poly s $A(x)$ and~$B(x)$ and let $n = A(1)B(1)$.  If
$$A(x)B(x) \equiv 1+x+\dots +x^{n-1} \pmod{x^n -1}$$
and  $p$ is a prime which is not a factor of~$A(1)$, then
$$A(x^p)B(x) \equiv 1+x+\dots +x^{n-1} \pmod{x^n -1}.$$
\endproclaim

\demo{Proof}  
Since $p$ is prime,
$A(x^p) \equiv \left(A(x)\right)^p \pmod p$, 
i.e., when the coefficients are reduced modulo~$p$. Let   
$G_n(x) = 1+x+\dots +x^{n-1}$. Then
$$A(x^p) B(x) = \left(A(x)\right)^{p-1}A(x)B(x) \equiv
  \left(A(x)\right)^{p-1}G_n(x),$$
where $\equiv$ means  the exponents are reduced modulo~$n$
and then the \coeff s  are reduced modulo~$p$.
Every $x^iG_n(x) \equiv G_n(x) \pmod{x^n-1}$, so
$$\left(A(x)\right)^{p-1} G_n(x) \equiv \left(A(1)\right)^{p-1}G_n(x)
  \pmod{x^n-1}.$$  
By Fermat's Little Theorem, 
$\left(A(1)\right)^{p-1} \equiv 1 \pmod{p}$.  Therefore 
$A(x^p)B(x) \equiv G_n(x)$, where the exponents are reduced modulo~$n$
and then the \coeff s are reduced modulo~$p$.  Both $A(x^p)B(x)$ and
$G_n(x)$ have \nneg\ coefficients whose sum is $n$ since
$A(1)B(1)=G_n(1)=n$.  Consider the following reductions.
\roster
 \item "{(R1)}" $A(x^p)B(x)$ is reduced modulo $x^n-1$, yielding a
        polynomial $G^*(x)$.
 \item "{(R2)}"The coefficients of $G^*(x)$ are reduced modulo $p$,
        yielding $G_n(x)$.
\endroster
(R1) preserves the sum of the coefficients, but (R2) reduces the sum
by some nonnegative multiple of $p$. Because the sum of the \coeffs\
of both $G^*(x)$ and $G_n(x)$ is~$n$, that multiple is $0$.  Therefore
$G^{*}(x)=G_n(x)$.
\qed \enddemo

We use the following result to prove Sands' Theorem (Lemma~2.4).  Let
$A-A$ be the difference set $\{a_1-a_2:a_1,a_2 \in A\}$.

\proclaim{Lemma 3.2} 
Let $A$ and $B$ be finite, $A,B \ne \{0\}$, and $A \oplus B$  a
complete set of residues modulo~$(\#A)(\#B)$.
Then at least one of the following  is true.
 \roster
   \item No member of $A-A$ is relatively prime to~$\#B$.
   \item No member of $B-B$ is relatively prime to~$\#A$.
 \endroster
\endproclaim

\demo{Proof}
Let $n = (\#A)(\#B)$. By Lemma~1.3,
$$A(x)B(x) \equiv 1 + x + \dots + x^{n-1}\pmod{x^n-1}.$$
Suppose $0 < a_1-a_2 = \delta'$ is relatively prime to~$\#B$ and 
$0 < b_1-b_2 =\delta''$ is
relatively prime to~$\#A$. Lemma~2.2 shows that
$$A(x^{\delta''})B(x^{\delta'}) \equiv  1 + x + \dots +
  x^{n-1}\pmod{x^n-1},$$
so by Lemma~1.3 again, $\delta'' A \oplus \delta' B$ is a complete set
of residues modulo~$n$. But
$$(b_1-b_2)a_1 + (a_1-a_2)b_2 =   (b_1-b_2)a_2 + (a_1-a_2)b_1.$$
Thus the same number  can be expressed $\delta'' a + \delta' b$ in two
ways, which is impossible.
\qed
\enddemo

\proclaim{Lemma 2.4} \cite{San} 
Let $A \oplus C=\Bbb Z$ be a tiling of period $n$ such that $A$ is
finite, $0 \in A \cap C$, and $n$ has one or two prime factors. Then
there is a prime factor $p$ of~$n$ such that either $A \subset p\Bbb Z$
or $C \subseteq p\Bbb Z$.
\endproclaim

\demo{Proof} 
Let $C=B \oplus n\Bbb Z$ and the prime factors of~$n$ be
$p$ and possibly~$q$.  Then at least one of 3.2(1) and~3.2(2) holds.

If 3.2(1) holds, then $A \subseteq A-A \subset p\Bbb Z \cup q \Bbb Z$,
the first containment because $0 \in A$. If neither $p\Bbb Z$
nor~$q\Bbb Z$ contains~$A$, then there exist $a_1,a_2 \in A$ such that
$a_1 \in p\Bbb Z \setminus q\Bbb Z$ and $a_2 \in q\Bbb Z \setminus
p\Bbb Z$.  But then $a_1-a_2$ is relatively prime to~$\#B$.

If 3.2(2) holds, the same argument shows that $B \subseteq p\Bbb Z$ or
$B \subseteq q\Bbb Z$. Then the same is true for $C=B \oplus n\Bbb Z$.
\qed \enddemo

\proclaim{Lemma 3.3}
Suppose $A$ is finite, $0 \in A$,  $A$ tiles the integers with
period~$n$, and $n$ has two prime factors, $p$ and~$q$.
If neither $\Phi_p(x)$ nor $\Phi_q(x)$ is a divisor of
$A(x)$, then $A \subset p\Bbb Z$ or $A \subset q\Bbb Z$.
\endproclaim

\demo{Proof}
Let $A \oplus (B \oplus n\Bbb Z)=\Bbb Z$
be a tiling of period~$n$.  By Lemma 1.3(4), $\Phi_p(x)$ and~
$\Phi_q(x)$ are divisors of $B(x)$.
From the remark after Lemma~1.4, 
neither $p \Bbb Z$ nor~$q \Bbb Z$
contains~$B$. Then the conclusion follows by Lemma~2.4.
\qed \enddemo

\head 4. A structure theory \endhead

In this section we describe the structure of those finite sets $A$
\st\ $A$ tiles the integers and $\#A$ has at most two prime factors.
Equivalently, such that the set $S_A$ of prime powers $s$ such that the
\cyclo\ \poly\ $\Phi_s(x)$ divides $A(x)$ consists of powers of
at most two primes. For $S$ such a set of prime powers, let $\Cal
T_S$ be the collection of all subsets $A$ of
$\{0,1,\dots,\lcm(S)-1\}$ which tile the integers and satisfy
$\min(A)=0$ and $S_A=S$. Note that $\Cal T_\varnothing = \{0\}$
because $\lcm(\varnothing )=1$, and  that $\Cal T_{\{p^{\alpha+1}\}}$
is the set whose only member is $p^{\alpha}\{0,1,\dots ,p-1\}$.
We have seen that there is no loss in requiring $\min (A)=0$.
We claim that a finite set $A'$ with $\min(A')=0$ and $S_{A'}=S$ tiles
the integers if and only if $A'$ is congruent modulo~$\lcm(S)$ to a
member of $\Cal T_S$. For if $A' \equiv A \pmod {\lcm(S)}$,
then $S_{A'} = S_{A} = S$, and as noted after the proof of Lemma~1.1,
$A' \oplus C_S=\Bbb Z$ if and only if $A \oplus C_S=\Bbb Z$.
Recall that $C_S$ is the universal translation set corresponding
to~$S$: $A \oplus C_S=\Bbb Z$ for every $A$ \st\ $A$ tiles the integers
and $S_A=S$.

For purposes of comparison we recall the simpler structure of {\it
all\/} finite sets which tile the \nneg\ integers $\Bbb N_0
=\{0,1,\dots\}$, due to deBruijn \cite{deB-3}. Note that every such set
has a unique translation set, so the unique associated tiling has a
period. One such set is $A=\{0,1,2,3,4\}\oplus \{0,10,20,30\} \oplus
\{0,120,240\}$, which tiles~$\Bbb N_0$ with period~$360$. $A$ can be
written $A=\tilde A \oplus 120\{0,1,2\}$, where $\tilde A$ tiles
~$\Bbb N_0$ with period~$60$, and it can be written
$A=\{0,1,2,3,4\} \oplus 5\bar A$, where $\bar A$ tiles~$\Bbb N_0$ with
period $72 =360/5$. If $A \ne \{0\}$ is any finite set which tiles~
$\Bbb N_0$, then there are always these two types of
direct sum decompostions, $A=\tilde A \oplus (n/p)\{0,1,\dots,p-1\}$
and $A=k\{0,1,\dots,q-1\} \oplus q\bar A$, where $p$ and~$q$
are prime factors of the period~$n$ of the tiling, $k=\gcd (A)$, and
$\tilde A$ and $\bar A$ are shorter tiles. Iterating either
decomposition, every tile is a direct sum, in one or more ways, of
tiles of the form $m\{0,1,\dots,p-1\}$. If the order is as above,
then $\tilde A$ is the direct sum of all but the last of the summands
and $q \bar A$ is the direct sum of all but the first. $A(x)$ is thus a
product of terms $(x^{mp}-1)/(x^m-1)=\Phi_p(x^m)$ and can easily be
shown to satisfy~(T1) and~(T2).

We return to  $\Cal T_S$ for the case that $S$ consists of the powers
of at most two primes. Both decompositions above generalize, the second
more usefully than the first.

Corresponding to the first decomposition, we will see that every member
of ~$\Cal T_S$ is a disjoint union of translates
of $(n/p)\{0,1,\dots,p-1\}$ and $(n/q)\{0,1,\dots,q-1\}$, where
$n=\lcm(S)$. The simplest case where both must be used is
$S=\{p,p^3,q^2\}$.  An important example of this with $\lcm(S)=72$ is
given below. More usefully, we will show that when 
$S \ne \varnothing$, every tile $A \in\Cal T_S$
is, as in Lemma~2.5, a union of translates of multiples of
$p$ or $q$ smaller tiles:
$A=m\bigcup_{i=0}^{p-1}\left(\{a_i\}\oplus p \bar A_i\right)$, where
$m=\gcd(A)$, $a_0=0$,
$\{a_0,a_1,\dots,a_{p-1}\}$ is a complete set of
residues modulo~$p$,
every  
$\{a_i\} \oplus p\bar A_i \subset \{0,1,\dots,\lcm(S)-1\}$,
and for some smaller set $\bar S$, every 
$\bar A_i\in \Cal T_{\bar S}$.
We need not get a direct sum, as the $\bar A_i$
need not be equal. Every $\bar A_i$ in turn is a union of
$p$ or~$q$ translates of multiples of even shorter tiles. Iterating the
procedure until $S=\varnothing$ gives the disjoint union referred to
above.

Suppose that  $S$ contains powers of only ~$p$, so that
$\lcm(S)$ is a power of $p$.
If $A \in \Cal T_S$, then $A \oplus C_S=\Bbb Z$ and either $p \in S$
and $C_S \subseteq p\Bbb Z$, or $p \notin S$ and $A \subset p\Bbb Z$.
Let $\bar S=\{p^{\alpha} : p^{\alpha +1} \in S\}$. If $p \notin S$,
then $\#\bar S=\#S$ and, as in Lemma~1.4, $\Cal T_S = \{p\bar A :
\bar A \in \Cal T_{\bar S}\}$.  If $p \in S$, then by the Corollary to
Lemma~2.5, the members of $\Cal T_S$ can be constructed by taking all
unions $\bigcup_{i=0}^{p-1}\left( \{a_i\} \oplus p\bar A_i\right)$
with $\bar A_i \in \Cal T_{\bar S}$,
$a_0=0$, $\{a_0,a_1,\dots,a_{p-1}\}$ a complete set of residues
modulo~$p$, and
every  
$\{a_i\} \oplus p\bar A_i \subset \{0,1,\dots,\lcm(S)-1\}$. 
This procedure gives all  of~$\Cal T_S$ and nothing else.

Suppose now that $S$ contains powers of both $p$ and $q$ and let
$$
 \bar S=\{p^{\alpha} : p^{\alpha +1} \in S\} \cup
  \{q^{\beta} : q^{\beta} \in S\}, \quad \bar S'=\{p^{\alpha} :
    p^{\alpha}\in S\} \cup \{q^{\beta} : q^{\beta+1} \in S\}.$$
We consider the three cases: $p \in S$, $q \in S$, and $p,q \notin S$. 
If $p \in S$, then $C_S \subseteq p\Bbb Z$ and $\Cal T_S$ can be
constructed as above by taking all unions 
$\bigcup_{i=0}^{p-1}\left(\{a_i\} \oplus p\bar A_i \right)$ 
with $\bar A_i \in \Cal T_{\bar S}$, $a_0 = 0$,
$\{a_0,a_1,\dots,a_{p-1}\}$ a complete
set of residues modulo~$p$, and
every 
$\{a_i\} \oplus p\bar A_i \subset \{0,1,\dots,\lcm(S)-1\}$.
If $q \in S$, then the analogous procedure, with the roles of
$p,\bar S$ and $q,{\bar S'}$ interchanged, gives $\Cal T_S$.  If both
$p$ and~$q$ are in~$S$, then $C_S \subseteq pq\Bbb Z$ and either
procedure gives $\Cal T_S$.  If neither~$p$ nor~$q$ is in~$S$, then
by Lemma~3.3, every member of~$\Cal T_S$ is contained in $p \Bbb Z$
or ~$q \Bbb Z$.  Then $\#S=\#\bar{S}=\#\bar S'$, and
$\{A \in \Cal T_S : A \subset p\Bbb Z\} = \{p\bar A : \bar A \in
  \Cal T_{\bar S}\}$,
while
$\{A \in \Cal T_S : A \subset q\Bbb Z\} =\{q\bar A : \bar A \in
  \Cal T_{\bar S'}\}$.
In all three cases, this procedure gives all of $\Cal T_S$ and nothing
else.

We  examine a few cases in more detail, including
the important example of deBruijn \cite{deB-2}.

When $S$ contains only powers of $p$, every member of $\Cal T_S$ is a
union of translates of $p^{\alpha}\{0,1,\dots,p-1\}$, where
$p^{\alpha+1}$ is the largest member of $S$. Hence every member of~
$\Cal T_S$ is a direct sum of this set and a set~$\tilde A$ which
also tiles the integers. Then
$S_{\tilde A}=S \setminus \{p^{\alpha+1}\}$,
but $\tilde A$ need not be a direct sum.
An example with $S =\{2,4,32\}$ is 
$16\{0,1\} \oplus \{0,1,2,11\}$.

It is easy to show that if $S = {\{p^{\alpha},q^{\beta}\}}$, then
every member of $ \Cal T_S$ is
$$p^{\alpha-1}\left(\tilde A \oplus
  pq^{\beta-1} \{0,1,\dots,q-1\}\right)$$
for $\tilde A \subseteq \{0,\dots,pq^{\beta-1}-~1\}$ 
a complete set of residues modulo~$p$ containing~$0$, or
an analogous set with the roles of $p$ and~$q$ interchanged.
Thus $\left\{k : \Phi_k(x) \text{ divides } A(x)\right\}$ contains
$\{p^\alpha\} \cup \{q^\beta,pq^\beta,p^2q^\beta,\dots,p^\alpha
 q^\beta\}$ or
$\{q^\beta\} \cup \{p^\alpha,p^\alpha q,p^\alpha q^2,\dots,p^\alpha
 q^\beta\}$.
If $\alpha>1$ and $\beta>1$, there are cyclotomic \poly\ divisors of
$A(x)$ in addition to the three required by~(T2).

The situation when $S$ has at least three elements is different. In
this case $\Cal T_S$ has members whose corresponding polynomial has
only the cyclotomic \poly\ divisors required by~(T2). We illustrate
this with the promised example. Among the members of 
$\Cal T_{\{4,9\}}$
are  $\bar A_0=\{0,3,6,18,21,24\}$ and
$\bar A_1 =\{0,2,12,14,24,26\} $. Each is a direct sum.
Consider
$$\align 
   A &=  \left(\{0\}\oplus 2\bar A_0 \right)
     \cup\left(\{1\} \oplus 2\bar A_1\right)\\
     &=   \{0,1,5,6,12,25,29,36,42,48,49,53\} \in \Cal T_{\{2,8,9\}}.
\endalign$$
The cyclotomic \poly\ divisors of $A(x)=\bar A_0(x^2)+x\bar A_1(x^2)$
are $\Phi_2(x)$ and those $\Phi_k(x)$ which divide both~ 
$\bar A_0(x^2)$ and~$\bar A_1(x^2)$, i\.e\.,
$$\align
   \left\{k : \Phi_k(x) \text{ divides } A(x)\right\} 
    &= \{2\} \cup \left(\{8,9,18,36,72\} \cap \{8,9,18,24,72\}\right)\\
    &= \{2,8,9,18,72\},
\endalign$$
exactly the set  required by~ (T2).
Then as in Theorem~A, $A \oplus (B \oplus 72\Bbb Z) = \Bbb Z$ for
$B=\{0,8,16,18,26,34\}$. deBruijn's example  was actually 
$\left(\{12\} \oplus 2\bar A_0\right)\cup\left(\{17\}
 \oplus 2\bar A_1\right)$.
It was the first example where
$A\oplus B$ is a complete set of  residues modulo~$n$
but neither $A$ nor~$B$ is periodic modulo~$n$.
Equivalently, 
neither $A$ nor~$B$ is
a disjoint union of translates of $(n/p)\{0,1,\dots,p-1\}$
for a single prime factor $p$ of~$n$.

\enddocument